\documentclass{article}
   \usepackage[OT4]{fontenc}
   \usepackage[english]{babel} 

   \usepackage[intlimits]{amsmath}
   \usepackage{amssymb}
   \usepackage{bbm}
   \usepackage{mathrsfs}
   \usepackage[nice]{nicefrac}

   \newcommand{\dx}{\,\text dx}
   \newcommand{\dy}[1]{\,\text d#1}
   \newcommand{\heat}[1]{\mathcal S(#1)}

   \def \calB{\mathcal B}

   \def \calT{{\mathcal T}}
   \def \calX{\mathcal X}

   \def \ee{{\,\text e}}
   \def \eps{\varepsilon}

   \def \R2{\mathbb R^2}

   \newtheorem{theorem}{Theorem}[section]
   \newtheorem{proposition}[theorem]{Proposition}
   \newtheorem{lemma}[theorem]{Lemma}
   \newtheorem{corollary}[theorem]{Corollary}

\begin{document}

   \begin{center}
   {\Large \textbf{Existence and asymptotics of solutions
   of the {D}ebye-{N}ernst-{P}lanck system in $\mathbb R^2$}}
   \vskip0.2in

      {  Agnieszka Herczak
         \footnote{Pa\'nstwowa Wy\.zsza Szko\a{l}a Zawodowa w Nysie, ul. Grodzka 19, 48-300 Nysa, Poland\\
            E-mail: {\tt herczak@math.uni.wroc.pl}},\quad
         Micha\a{l} Olech
         \footnote{Instytut Matematyczny, Uniwersytet Wroc\a{l}awski, pl. Grunwaldzki 2/4, 50-384 Wroc\a{l}aw, Poland;
            E-mail: {\tt olech@math.uni.wroc.pl},\\
            The preparation of this article has been partially supported by a KBN/MNiI grant \texttt{1 P03A 008 30}}}
\end{center}

   \abstract{In this paper we investigate a system describing electrically charged particles in the whole space $\mathbb R^2$. Our main goal is to describe large time behavior of solutions which start their evolution from initial data of small size. This is achieved using radially symmetric self-similar solutions.}

\vskip 10pt

\hskip-\parindent \textbf{Key words:}\ parabolic-elliptic system, existence of solutions, asymptotic behavior of solutions, long time behavior of solutions.

\vskip 10pt

\hskip-\parindent \textbf{AMS subject classification:}\ Primary 35B40, 35K15, 35Q99; Secondary 78A35, 92E99.

\vskip 10pt

\hskip-\parindent \textbf{Acknowledgment:}\ Authors are grateful to Professors\ Piotr\ Biler and Grzegorz\ Karch for a very fruitful discussions and their invaluable help.

\section{Introduction}

We analyze the Debye-Nernst-Planck system of the form 
\begin{subequations}\label{AGHproblem}
   \begin{equation}
   \begin{split}\label{AGHmainintro}
      u_t &=\Delta u + \nabla \cdot \big(u\nabla \phi_{u}\big), \\
      \phi_{u} &= - E_2\ast u\text{\quad in \quad}\mathbb R^2,
   \end{split}
   \end{equation}
where $u(x,t): \mathbb R^2\times \mathbb R^+ \rightarrow \mathbb R$ and $E_2(z)=\tfrac{1}{2\pi}\log|z|$ is the fundamental solution of the Laplace equation in $\mathbb R^2$. We supplement this system of equations with the initial condition
   \begin{equation}\label{AGHinitintro}
      u(x,0) = u_0(x)  \text{\quad in\quad}\mathbb R^2.
   \end{equation}
\end{subequations}
One may think about this system as of a single parabolic equation with a nonlocal nonlinearity.

The history of studies of this model is quite long. The system was introduced by Nernst and Planck in the nineteenth century, and then modified by Debye in twenties of the twentieth century (see \cite{debye_huckel}). It describes the evolution of the density $u(x,t)$ of electrically charged particles interacting through the potential $\phi_u$. We encounter systems of this type, also for several unknown densities,  in the electrochemistry, theory of semiconductors and physics of plasma \cite{babskii,choi_lui,jungel,markowich_schmeiser}.

Many different methods have been used for the analysis of this system since that time (see \cite{biler_92, biler_95, biler_dolbeault, biler_hebisch_nadzieja} and references therein). However some problems still remain unsolved.

Similar problems where the first equation of \eqref{AGHmainintro} is replaced by $\phi_u=E_2\ast u$, appear in the chemotaxis theory and they are used for describing the evolution of particles which interact via gravitational forces \cite{biler_nadzieja,horstmann_1, horstmann_2}. Some reasonings from that theory can be adopted without essential changes but questions related to the long time asymptotics have completely different answers. For example, Banach type theorems giving the local existence can be applied in both cases (see the next section). On the other hand,  we know that gravitational problem does not have solutions for large initial masses (see \cite{suzuki_book} and the references therein), which is not the case of electric forces.

Whenever it would not cause any confusion we omit time or spatial variable. Moreover $(f\ast g)\,(x)$ denotes the usual convolution, i.e.\ $(f\ast g)\,(x) = \int_{\mathbb R^2} f(x-y)g(y)\dy{y}$, and  $\heat{t}$ is the heat semigroup given by
   \begin{equation}\label{AGHpolgrupa}
      \big(\heat{t}u \big)(x)=\int_{\mathbb R^2} G(x-y,t)u(y)\dy{y},
   \end{equation}
where $G(x,t)=\tfrac{1}{4\pi t}\,\ee^{-\nicefrac{|x|^2}{4t}}$ is the Gauss kernel. Here and subsequently $\|u\|_p$ denotes the usual $L^p(\mathbb R^2)$ norm of the function $u(x,t)$, with respect to the spatial variable. $C$ is a generic constant which can vary from line to line. The Schwartz class of rapidly decreasing functions is denoted by $\mathscr S(\mathbb R^2)$ and its dual by $\mathscr S'(\mathbb R^2)$. $\delta_0(x)$ denotes Dirac delta function at  the point $0$.\medskip

In each dimension if $u(x,t)$ and $\phi(x,t)$ solve the system \eqref{AGHmainintro} then $\lambda^2u(\lambda x,\lambda^2 t)$ and $\phi(\lambda x,\lambda^2 t)$ also do. As we will see the total charge of particles, namely the integral $\int_{\mathbb R^2} u(x,t)\dx$, is preserved in time. However, solutions invariant under the above scaling preserve the charge only in $\mathbb R^2$. So that we expect essentially nonlinear asymptotics  of generic solutions. Moreover, self-similar solutions which  describe long time behavior decay faster than the Gauss kernel when $|x|\rightarrow\infty$. The two dimensional case differs much from the higher dimensional case, where the long time behavior is described simply by the heat kernel \cite{biler_dolbeault}. However, in \cite {biler_dolbeault} authors use a completely different approach involving energy and entropy dissipation methods.

First we prove (following methods used in \cite{biler_95,karch}) the existence of solutions for the system \eqref{AGHmainintro}. Then we study self-similar solutions in radially symmetric and general setting. The last section is devoted to the large time behavior of solutions.

\section{Existence of solutions for singular initial data}\label{AGHmildSol}

We are mainly interested in {\it mild} solutions framework. We look for solutions of the integral equation 
   \begin{equation}\label{AGHcalkowe}
      u(t)=\heat{t}u_0+\int_0^t \heat{t-\tau}\, \nabla \cdot \big(u \nabla \phi \big)(\tau)\dy{\tau},
   \end{equation}
where $\phi(x,t)$ and $u_0(x)$ are as in \eqref{AGHproblem}. Two major problems connected with this approach appear. One is: in what sense the solutions of the integral equation \eqref{AGHcalkowe} satisfy the differential problem \eqref{AGHmainintro}. It is obvious that the classical solutions of \eqref{AGHmainintro} satisfy the equality \eqref{AGHcalkowe}. Unfortunately, the reverse implication is not always true (see also a discussion of this problem in \cite{karch}).\\
The second problem is: in what sense the initial condition is fulfilled. The answer for this question will come with the definition of the space in which we solve the integral problem \eqref{AGHcalkowe}.

Let us notice here that the term $\nabla\cdot \big( u \nabla \phi_u )$ can be written as $\calB(u,u) = \nabla \cdot \big( u \nabla \phi_u \big)$ where $\calB(u,v)$ is a bilinear form defined by 
\begin{equation*}
   \calB(u,v) = \nabla \cdot \big( u \nabla \phi_v \big),\quad \text{with} \quad \phi_v = - E_2 \ast v,
\end{equation*}
for $u,v$ from some Banach space $\mathcal X$, which in our case is a space $L^{\nicefrac{4}{3}}(\R2)$ (see a discussion before a proof of Lemma \ref{AGHapriori} below).\\
We introduce here the notation of \textit{the rescaled version} of a given function $f(x)$
\begin{equation*}
   f_\lambda(x) = f(\lambda x),
\end{equation*}
where $\lambda$ is some real and positive parameter. In the case of the function depending on space and time, this kind of scaling (and, therefore, the notation) is used only with respect to the spatial variable.

By a simple calculation we show that the bilinear form $\calB(u,v)$ has scaling order equal to 0 (cf. \cite{karch}) which means that for every $\lambda > 0$ and every $u,v \in \calX$
\begin{equation*}
   \calB \big( u_\lambda,v_\lambda \big) = \big( \calB (u,v) \big)_\lambda.
\end{equation*}

We use frequently very well known estimates of the heat semigroup 
\begin{equation}\label{AGHinheat}
\begin{split}
   &\big \|\heat{t}f \big \|_{p}\leqslant C t^{\nicefrac{1}{p} - \nicefrac{1}{q}} \|f\|_{q}\hskip 5pt,\\
   &\big \|\nabla \heat{t}f \big \|_p \leqslant C t^{\nicefrac{1}{p} - \nicefrac{1}{q} - \nicefrac{1}{2}}\|f\|_q \hskip 5pt,
\end{split}
\end{equation}
where $1\leqslant q \leqslant p \leqslant\infty$, $C=C(p,q)$ is some positive constant, $t>0$ and $f(x)\in L^p(\R2)$. Let us also recall the weak Young inequality 
   \begin{equation}\label{AGHweaky}
      \|g\ast h\|_{\gamma}\leqslant C({\beta},{\gamma},g)\|h\|_{\alpha}\hskip 3pt.
   \end{equation}
for all $h\in L^{\alpha}(\mathbb R^2)$, $|g(x)|\leqslant |x|^{-\nicefrac{2}{\beta}}$ and $\nicefrac{1}{\alpha}+\nicefrac{1}{\beta} = 1 + \nicefrac{1}{\gamma}$. For the proof see, e.g., \cite[Remark (2) in Theorem 4.3]{lieb_loss}. 

The inclusions $\mathscr S(\mathbb R^2) \subset L^{\nicefrac{4}{3}}(\R2) \subset \mathscr S'(\mathbb R^2)$ are continuous and the norm $\|\cdot \|_{\nicefrac{4}{3}}$ is translation invariant. Moreover the estimation
\begin{equation}\label{AGHszacB}
   \big \| \heat{t}\calB(u,v) \big \|_{\nicefrac{4}{3}} \leqslant C\,t^{-\nicefrac{3}{4}} \|u\|_{\nicefrac{4}{3}}\|v\|_{\nicefrac{4}{3}}
\end{equation}
holds, for every $u,v \in L^{\nicefrac{4}{3}}(\R2)$.\\
Indeed, using the second inequality in \eqref{AGHinheat} with $p = \nicefrac{4}{3}$ and $q = 1$ and the H\"older inequality, we get
\begin{equation*}
   \big \| \heat{t}\calB(u,v) \big \|_{\nicefrac{4}{3}} \leqslant C\,t^{-\nicefrac{3}{4}} \|u\|_{\nicefrac{4}{3}} \|\nabla \phi_v\|_{4}.
\end{equation*}
We use then inequality \eqref{AGHweaky} with $\gamma = 4$, $\beta = 2$ and $\alpha = \nicefrac{4}{3}$ to obtain
\begin{equation*}
   \|u\|_{\nicefrac{4}{3}} \|\nabla \phi_v\|_{4} = \|u\|_{\nicefrac{4}{3}} \big \| |x|^{-1} \ast v \big \|_{4} \leqslant \|u\|_{\nicefrac{4}{3}} \|v\|_{\nicefrac{4}{3}},
\end{equation*}
which concludes the proof of \eqref{AGHszacB}.

\medskip

It means that the space $L^{\nicefrac{4}{3}}(\R2)$ is {\it adequate} to the problem \eqref{AGHproblem} in the sense of \cite[Definition 3.1]{karch}.
We construct a homogeneous Besov space modeled on the Banach space $\mathcal X = L^{\nicefrac{4}{3}}(\R2)$, as
\begin{equation*}
   B\calX = \Big \{ u \in\mathscr S'(\R2):\ \sup_{t > 0} t^{\nicefrac{1}{4}} \big \|\heat{t} u \big \|_{\nicefrac{4}{3}} < \infty \Big \}\hskip 5pt,
\end{equation*}
and we introduce a norm on this space by
\begin{equation*}
   \| u \|_{B\calX} = \sup_{t > 0} t^{\nicefrac{1}{4}} \big \|\heat{t} u \big \|_{\nicefrac{4}{3}}\hskip 5pt .
\end{equation*}
The space $B\calX$ is in fact the classical homogeneous Besov space $\dot{B}^{\nicefrac{1}{2}}_{\nicefrac{4}{3},\infty}$. The norm $\|\cdot\|_{B\calX}$ is equivalent to the original norm in $\dot{B}^{\nicefrac{1}{2}}_{\nicefrac{4}{3},\infty}$ introduced by the dyadic decomposition. For a further discussion on this topic we refer the reader to \cite[Example 4.2]{karch} and the references therein.

As it was in {\cite[Lemma 5.2]{karch}} let us prove here the \textit{a~priori} estimate for the bilinear form $\calB(u,v)$.
\begin{lemma}\label{AGHapriori}
   There exists a constant $C > 0$ independent of \hskip 3pt $u,\,v$ and $t$ such that 
   \begin{equation}\label{AGHapreq}
      \big \| \heat{t}\calB(u,v) \big \|_{B\calX} \leqslant C\,t^{-\nicefrac{1}{2}} \|u\|_{\nicefrac{4}{3}} \, \|v\|_{\nicefrac{4}{3}}\hskip 5pt,
   \end{equation}
   for all $u,v \in \calX$.
\end{lemma}

{\hskip -\parindent \sc Proof\hskip10pt}
   We begin with the fact, that there exists some function $C(t)$, independent of $u$ and $v$ such that
\begin{equation*}
      \big \| \heat{t}\calB(u,v) \big \|_{B\calX} = \sup_{s>0} s^{\nicefrac{1}{4}} \big \| \heat{s}\heat{t}\calB(u,v) \big \|_{\nicefrac{4}{3}} \leqslant C(t) \|u\|_{\nicefrac{4}{3}} \, \|v\|_{\nicefrac{4}{3}}
   \end{equation*}
    for every $t > 0$. The proof consists of two steps. First let us assume that $0 < s < 1$. Then also $0 < s^{\nicefrac{1}{4}} < 1$ and because $\|\heat{t}f\|_{\nicefrac{4}{3}} \leqslant \|f\|_{\nicefrac{4}{3}}$, the inequality \eqref{AGHszacB} yields
   \begin{equation*}
      s^{\nicefrac{1}{4}} \big \| \heat{s}\heat{t}\calB(u,v) \big \|_{\nicefrac{4}{3}} \leqslant \big \| \heat{t}\calB(u,v) \big \|_{\nicefrac{4}{3}} \leqslant C t^{-\nicefrac{3}{4}}\| u \|_{\nicefrac{4}{3}}\,\| v \|_{\nicefrac{4}{3}}\hskip 5pt.
   \end{equation*}
   In the case $s > 1$ we proceed in a similar way but first we use the semigroup property $\heat{s}\heat{t} = \heat{s+t} = \heat{t}\heat{s}$ to obtain
   \begin{gather*}
      s^{\nicefrac{1}{4}} \big \| \heat{s}\heat{t}\calB(u,v) \big \|_{\nicefrac{4}{3}} = s^{\nicefrac{1}{4}} \big \| \heat{t} \heat{s} \calB(u,v) \big \|_{\nicefrac{4}{3}} \\
      \leqslant s^{\nicefrac{1}{4}} \big \| \heat{s} \calB(u,v)\big \|_{\nicefrac{4}{3}} \leqslant C s^{-\nicefrac{1}{2}}\| u \|_{\nicefrac{4}{3}}\,\| v \|_{\nicefrac{4}{3}} \leqslant C \| u \|_{\nicefrac{4}{3}}\,\| v \|_{\nicefrac{4}{3}}\hskip 5pt.
   \end{gather*}
   Now, we show that $C(t)$ is some positive constant. First let us check the scaling property
   \begin{equation}\label{AGHscalStB}
      \heat{t}\calB(u_\lambda,v_\lambda) = \big(\heat{\lambda^2 t} \calB(u,v) \big)_\lambda \hskip 5pt,
   \end{equation}
   for all $u,\ v\in L^{\nicefrac{4}{3}}(\R2)$ and for all real $\lambda >0$.
   For the proof let $\xi\in\mathscr S'(\R2)$. Then
   \begin{equation}\label{AGHscheatb}
   \begin{split}
      \big( \heat{t}\xi_\lambda \big)(x) &= \frac{1}{2\pi t}\int_{\R2} \ee^{-\frac{|x-z|^2}{4t}}\xi(\lambda z)\dy{z} = \frac{1}{2\pi \lambda^2 t}\int_{\R2} \ee^{-\frac{\lambda^2|x-y/\lambda|^2}{4\lambda^2 t}}\xi(y)\dy{y}\\
      &= \frac{1}{2\pi \lambda^2 t}\int_{\R2} \ee^{- \frac{|\lambda x - y|^2}{4\lambda^2 t}}\xi(y)\dy{y} = \big( \heat{\lambda^2 t}\xi(x) \big)_\lambda \hskip 5pt .
   \end{split}
   \end{equation}
   As we mentioned before the scaling order of the form $\calB(u,v)$ equals zero which yields
   \begin{align*}
      \heat{t} \calB(u_\lambda,v_\lambda) = \heat{t}\big( \calB(u,v) \big)_\lambda = \Big( \heat{\lambda^2 t}\calB (u,v)\Big)_\lambda.
   \end{align*}
   Then we use scaling properties
   \begin{equation*}
      \|f_\lambda\|_{\nicefrac{4}{3}} = \lambda^{-\nicefrac{3}{2}} \|f\|_{\nicefrac{4}{3}} \quad \text{and} \quad \|f\|_{B\calX} =  \lambda^{2}\|f_\lambda \|_{B\calX} \hskip 3pt,
   \end{equation*}
   to show that
   \begin{gather*}
      \big \| \heat{\lambda^2 t}\calB(u,v) \big \|_{B\calX} = \lambda^{2} \big \| \big (\heat{\lambda^2 t} \calB(u,v) \big )_\lambda \big \|_{B\calX} = \lambda^{2} \big \| \heat{t} \calB(u_\lambda,v_\lambda) \big \|_{B\calX} \\
      \leqslant \lambda^{2} C(t) \|u_\lambda \|_{\nicefrac{4}{3}} \, \|v_\lambda \|_{\nicefrac{4}{3}} = C(t) \lambda^{-1}\| u \|_{\nicefrac{4}{3}} \| v \|_{\nicefrac{4}{3}}\hskip 5pt.
   \end{gather*}
   Now we fix $t = t_0 > 0$ and we put $\lambda = (\nicefrac{t}{t_0})^{\nicefrac{1}{2}}$. Then we obtain the estimate
   \begin{equation*}
      \big \| \heat{t}\calB(u,v) \big \|_{B\calX} \leqslant C(t_0)\,t_0^{\nicefrac{1}{2}} t^{-\nicefrac{1}{2}} \| u \|_{\nicefrac{4}{3}} \| v \|_{\nicefrac{4}{3}},
   \end{equation*}
   which concludes the proof. \hskip 10pt $\blacksquare$ \medskip

The space where we are looking for the solutions of the integral equation \eqref{AGHcalkowe} is defined as
\begin{equation}\label{AGHspaceX}
\begin{split}
   \mathscr X = \mathscr F \big( [0,\infty): & B\calX \big)\ \cap\\
   &\cap\, \Big \{ f(t):\ (0,\infty) \rightarrow L^{\nicefrac{4}{3}}(\R2):\ \sup_{t > 0} t^{\nicefrac{3}{4}} \|f(t)\|_{\calX} < \infty \Big \},
\end{split}
\end{equation}
and we equip this space with the norm
\begin{equation*}
   \| f \|_{\mathscr X} = \max \Big \{ \sup_{t >0} \|f(t)\|_{B\calX},\ \sup_{t >0} t^{\nicefrac{3}{4}}\|f(t)\|_{\nicefrac{4}{3}} \Big \}.
\end{equation*}
The space $\mathscr X$ with the above norm is a Banach space. We understand $\mathscr F \big( [0,\infty): B\calX \big)$ as the space of $B\calX$-valued measurable functions which belong to $L^\infty\big( [0,\infty): B\calX \big)$ and attain the initial condition of the problem \eqref{AGHproblem} as $t \rightarrow 0^+$ in the sense of tempered distributions, namely
\begin{equation*}
   \int_{\R2} \phi(x) u(x,t) \dx \rightarrow \int_{\R2} \phi(x)u_0(x) \dx,
\end{equation*}
for all $\phi \in \mathscr S(\mathbb R^2)$ as $t$ goes to $0^+$.

Now let us formulate and prove the main result of this section. In the proof we combine the reasonings presented in \cite{biler_95} and \cite{karch, karch2}.

\begin{theorem}\label{AGHGlobalExistence}
   There exists an $\eps > 0$ such that for all $u_0 \in B\calX$ satisfying $\|u_0\|_{B\calX} < \eps$, there is a global in time solution $u(x,t)$ of the problem \eqref{AGHproblem} in the space \eqref{AGHspaceX}.
   This solution is unique among all the functions $u\in B\calX$ satisfying the condition $\sup_{t > 0} t^{\nicefrac{1}{4}} \|u(t)\|_{\calX} < 2\eps$.
\end{theorem}

{\hskip -\parindent \sc Proof \hskip 7pt}
   We show that the operator
   \begin{equation*}
      \mathscr N \big(u \big)(t) = \heat{t}u_0 + \int_0^t \heat{t-\tau}\, \calB \big( u,u \big) (\tau) \dy{\tau}
   \end{equation*}
   is a contraction in the set
   \begin{equation*}
      B_{\eps} = \Big \{ u \in \mathscr X:\ \|u(t)\|_{\mathscr X} \leqslant 2 \eps \Big \}\hskip 5pt
   \end{equation*}
   for some small $\eps >0$, to be specified later.\\
   In the following estimates we use frequently the elementary inequalities such as the H\"older inequality, the heat semigroup estimates \eqref{AGHinheat} and the estimates of the bilinear form $\calB(u,v)$ from \eqref{AGHszacB} and \eqref{AGHapreq}
   \begin{align*}
      t^{\nicefrac{1}{4}} \big \| \mathscr N \big( u \big) (t) \big \|_{\nicefrac{4}{3}} \leqslant&\, \big \| u_0 \big \|_{B\calX} + C\,t^{\nicefrac{1}{4}} \int_0^t (t-\tau)^{-\nicefrac{3}{4}} \tau^{-\nicefrac{1}{2}} \bigg( \tau^{\nicefrac{1}{4}} \big \|u(\tau)\big \|_{\nicefrac{4}{3}} \bigg)^2  \dy{\tau} \\
      \leqslant&\, \big \| u_0 \big \|_{B\calX} + C\,\bigg( \sup_{\tau >0} \tau^{\nicefrac{1}{4}} \big \|u(\tau) \big \|_{\nicefrac{4}{3}} \bigg)^2 \int_0^1 (1 - w)^{-\nicefrac{3}{4}} w^{-\nicefrac{1}{2}} \dy{w} \\
      \leqslant&\, \big \| u_0 \big \|_{B\calX} + C \eps^2,
   \end{align*}
   and similarly
   \begin{multline*}
      t^{\nicefrac{1}{4}} \big \| \mathscr N \big( u \big) (t) - \mathscr N \big( v \big) (t) \big \|_{\nicefrac{4}{3}}\\
      \leqslant\ C\,t^{\nicefrac{1}{4}} \int_0^t (t-\tau)^{-\nicefrac{3}{4}} \big \| (u-v)(\tau) \big \|_{\nicefrac{4}{3}} \Big( \big \| u(\tau) \big \|_{\nicefrac{4}{3}} + \big \| v(\tau) \big \|_{\nicefrac{4}{3}}\Big) \dy{\tau}\\
      \leqslant\ C\, \eps \sup_{t > 0} t^{\nicefrac{1}{4}} \big \| (u-v)(t) \big \|_{\nicefrac{4}{3}}\hskip 5pt.
   \end{multline*}
   In the same way we deal with the estimate
   \begin{align*}
      \big \| \mathscr N \big( u \big) (t) \big \|_{B\calX} \leqslant&\, \big \| u_0 \big \|_{B\calX} + C \int_0^t (t-\tau)^{-\nicefrac{1}{2}} \tau^{-\nicefrac{1}{2}} \bigg( \tau^{\nicefrac{1}{4}} \big \|u(\tau) \big \|_{\nicefrac{4}{3}} \bigg)^2  \dy{\tau}\\
      \leqslant&\, \big \| u_0 \big \|_{B\calX} + C \bigg( \sup_{\tau > 0} \tau^{\nicefrac{1}{4}} \big \|u(\tau) \big \|_{\nicefrac{4}{3}} \bigg)^2 \int_0^1 (1 - w)^{-\nicefrac{1}{2}} w^{-\nicefrac{1}{2}} \dy{w}\\
      \leqslant&\, \big \| u_0 \big \|_{B\calX} + C \eps^2
   \end{align*}
   and
   \begin{equation*}
      \big \| \mathscr N \big( u \big) (t) - \mathscr N \big( v \big) (t) \big \|_{B\calX} \leqslant \, C\, \eps \sup_{t > 0} t^{\nicefrac{1}{4}} \big \| (u-v)(t) \big \|_{\nicefrac{4}{3}}\hskip 5pt.
   \end{equation*}
   Combining the above inequalities and choosing $\eps >0$ such that $C \eps^2 < \eps$ and $C \eps < 1$, we obtain
   \begin{gather*}
      \big \| \mathscr N (u) \big \|_{\mathscr X} \leqslant 2 \eps \hskip 5pt,\\
      \big \| \mathscr N (u) - \mathscr N (v) \big \|_{\mathscr X} \leqslant \big \| (u-v)(t) \big \|_{\mathscr X} \hskip 5pt,
   \end{gather*}
   for every  $u,\ v\in L^{\nicefrac{4}{3}}(\R2)$. It follows that $\mathscr N:\ B_{\eps} \rightarrow B_{\eps}$ is a contraction in $B_{\eps}$. Thanks to the Banach fixed point theorem we obtain the existence of a mild solution of the problem \eqref{AGHproblem}. This solution is unique in the set $B_{\eps}$ and is global in time. {\hskip 10pt $\blacksquare$} \medskip

Since we are particularly interested in self-similar solutions, it is desired to know what does the condition $\|u_0\|_{B\calX} < \eps$ mean for the initial condition being a distribution. Namely, let $u_0(x)=M\delta_{x_0}(x)$. Then
\begin{equation*}
\begin{split}
   \|M\heat{t}\delta_{x_0}\|_{\nicefrac{4}{3}}^{\nicefrac{4}{3}}=&\ \int_{\mathbb R^2} \Bigg( \frac{M}{4\pi t}\int_{\mathbb R^2} \ee^{-\frac{|x-y|^2}{4t}}\,\delta_{x_0}(y)\dy{y}\Bigg)^{\nicefrac{4}{3}}\dx\\
      =&\ \int_{\mathbb R^2} \bigg(\frac{M}{4\pi t}\;\ee^{-\frac{|x-x_0|^2}{4t}} \bigg)^{\nicefrac{4}{3}}\!\!\!\dx = \frac{M^{\nicefrac{4}{3}}}{(4\pi t)^{\nicefrac{1}{3}}} \int_{\mathbb R^2} \frac{1}{4\pi t}\; \ee^{-\frac{4|x-x_0|^2}{3 \cdot 4t}}\dx\\
      =&\ \frac{4}{3} \cdot \frac{M^{\nicefrac{4}{3}}}{(4\pi t)^{\nicefrac{1}{3}}} \cdot \frac{1}{t^{\nicefrac{1}{4}}} \int_{\mathbb R^2} \frac{1}{4\pi t}\; \ee^{-\nicefrac{|z|^2}{4t}}\dy{z}.
\end{split}
\end{equation*}
Then
   \begin{equation*}
      \limsup_{t \rightarrow 0} t^{\nicefrac{1}{4}} \|M \heat{t}\delta_{x_0}\|_{\nicefrac{4}{3}} \sim M^{\nicefrac{4}{3}},
   \end{equation*}
which means that atomic measures are admissible as initial data \eqref{AGHinitintro} but only with small masses of single atoms.

Moreover let us notice, that if $u_0(x)$ is a function from $L^1(\mathbb R^2)$ -- it does not contain singular part -- then
\begin{equation*}
   \limsup_{t \rightarrow 0} t^{\nicefrac{1}{4}} \| u_0 \|_{\nicefrac{4}{3}} = 0.
\end{equation*}
Above shows that the part $\limsup_{t \rightarrow 0} t^{\nicefrac{1}{4}} \| \cdot \|_{\nicefrac{4}{3}}$ of the norm $\|\cdot\|_{\mathscr X}$ measures only the singular part.

\section{Local and global existence in $L^p(\R2)$ space} \label{AGHmildLp}

Using the same framework of mild solutions we prove the local in time existence of solutions of the problem \eqref{AGHcalkowe} in the space $C([0,T];L^1(\R2)) \cap C([0,T];L^p(\R2))$ for some $p>1$ and an arbitrary initial condition $u_0 \in L^1(\R2)\cap L^p(\R2)$. Moreover, we show that these solutions can be extended to global in time solutions. This situation differs from the case where particles interact through gravitational forces and through the chemotactic attraction. The solutions in that case blow up in $L^p(\R2)$ spaces for sufficiently large initial data (see \cite{horstmann_1, horstmann_2}).

In the following we use the same inequalities as before, i.e.\ the heat semigroup estimates \eqref{AGHinheat} and the weak Young inequality \eqref{AGHweaky}.

Let us define the space
\begin{equation*}
   \calX_T = C \big([0,T];L^1(\R2) \big) \cap C \big([0,T];L^p(\R2) \big),
\end{equation*}
which is Banach space with the norm
\begin{equation*}
   \|u\|_{\calX_T} = \sup_{0\leqslant t \leqslant T} \|u\|_1 + \sup_{0\leqslant t \leqslant T} \|u\|_p \hskip 5pt.
\end{equation*}

\begin{proposition}\label{AGHlocex}
   Let $\frac{4}{3} \leqslant p < 4$. For any $u_0 \in L^1(\R2)\cap L^p(\R2)$ there exists some $T>0$ and a unique solution of the integral equation
   \begin{equation*}
      u(t) = \heat{t} u_0 + \int_0^t \heat{t - \tau}\nabla \cdot \big(u \nabla \phi \big)(\tau) \dy{\tau},
   \end{equation*}
   such that $u(x,t) \in \calX_T$. 
\end{proposition}

{\hskip -\parindent \sc Proof \hskip 10pt}
   Let 
   \begin{equation}\label{AGHOperLoc}
      \mathscr N \big( u \big) (t) = \heat{t}u_0 + \int_0^t \nabla \heat{t-\tau}\cdot \big(u \nabla \phi \big)(\tau) \dy{\tau}.
   \end{equation}
   Then \eqref{AGHinheat} implies that
   \begin{equation*}
      \mathcal I_p(t) = \|\mathscr N \big( u \big) (t) - \heat{t}u_0\|_p \leqslant C \int_0^t (t - \tau)^{\nicefrac{1}{p} - \nicefrac{1}{q} -\nicefrac{1}{2}} \big \| \big(u \nabla \phi \big)(\tau) \big \|_q \dy{\tau},
   \end{equation*}
   for every $1\leqslant p,\, q \leqslant \infty$. We use the H{\"o}lder inequality to obtain
   \begin{equation*}
      \mathcal I_p(t) \leqslant C \int_0^t (t - \tau)^{\nicefrac{1}{p} - \nicefrac{1}{q} -\nicefrac{1}{2}} \big \| u(\tau) \big \|_p \big \| \big( \nabla \phi \big)(\tau) \big \|_{\frac{pq}{p-q}} \dy{\tau}.
   \end{equation*}
   We can exploit \eqref{AGHweaky} with $\alpha = \frac{2pq}{pq - 2p - 2q}, \beta = 2$ and $\gamma = \frac{pq}{p-q}$, which yields
   \begin{equation*}
      \mathcal I_p(t) \leqslant C \int_0^t (t - \tau)^{\nicefrac{1}{p} - \nicefrac{1}{q} -\nicefrac{1}{2}} \big \| u(\tau) \big \|_p \big \| u(\tau) \big \|_{\frac{2pq}{pq - 2p - 2q}} \dy{\tau}.
   \end{equation*}
   We would like to have $\frac{2pq}{pq - 2p - 2q} = p$ which is possible for $q= \frac{2p}{4-p}$ and all $\nicefrac{4}{3} \leqslant p < 4$. Consequently, we have
   \begin{multline*}
      \sup_{0 \leqslant t \leqslant T} \|\mathscr N \big( u \big) (t) - \heat{t}u_0\|_p \\
      \leqslant C \bigg( \sup_{0 \leqslant t \leqslant T} \big \| u(t) \big \|_p \bigg)^2 \int_0^t (t - \tau)^{-\nicefrac{1}{p}} \dy{\tau}\\
      \leqslant C T^{1-\nicefrac{1}{p}} \bigg( \sup_{0 \leqslant t \leqslant T} \big \| u(t) \big \|_p \bigg)^2 \int_0^1 (1 - s)^{-\nicefrac{1}{p}} \dy{s},
   \end{multline*}
   where the last integral is finite since $- \nicefrac{3}{4} \leqslant - \nicefrac{1}{p} \leqslant -\nicefrac{1}{4}$.

   Similarly, we obtain
   \begin{multline*}
      \sup_{0 \leqslant t \leqslant T} \|\mathscr N \big( u \big) (t) - \heat{t}u_0\|_1 \\
      \leqslant C \bigg( \sup_{0 \leqslant t \leqslant T} \big \| u(t) \big \|_p \bigg)^2 \int_0^t (t - \tau)^{1 - \nicefrac{2}{p}} \dy{\tau}\\
      \leqslant C T^{2 - \nicefrac{2}{p}} \bigg( \sup_{0 \leqslant t \leqslant T} \big \| u(t) \big \|_p \bigg)^2 \int_0^1 (1 - s)^{1- \nicefrac{2}{p}} \dy{s},
   \end{multline*}
   with the same constraint for $p$, namely $\nicefrac{4}{3} \leqslant p < 4$, which yields the inequality $- \nicefrac{1}{2} \leqslant 1 - \nicefrac{2}{p} \leqslant \nicefrac{1}{2}$. \\
   An analogous calculation gives
   \begin{multline*}
      \sup_{0 \leqslant t \leqslant T} \|\mathscr N \big( u \big) (t) - \mathscr N \big( v \big) (t)\|_p \\
      \leqslant C T^{1-\nicefrac{1}{p}} \sup_{0 \leqslant t \leqslant T} \big \| u(t) - v(t)\big \|_p \bigg(\sup_{0 \leqslant t \leqslant T} \big \| u(t) \big \|_p + \sup_{0 \leqslant t \leqslant T} \big \| v(t) \big \|_p \bigg)
   \end{multline*}
   and
   \begin{multline*}
      \sup_{0 \leqslant t \leqslant T} \|\mathscr N \big( u \big) (t) - \mathscr N \big( v \big) (t)\|_1 \\
      \leqslant C T^{2-\nicefrac{2}{p}} \sup_{0 \leqslant t \leqslant T} \big \| u(t) - v(t)\big \|_p \bigg(\sup_{0 \leqslant t \leqslant T} \big \| u(t) \big \|_p + \sup_{0 \leqslant t \leqslant T} \big \| v(t) \big \|_p \bigg),
   \end{multline*}
   for all $u,v \in \calX_T$ and $\nicefrac{4}{3} \leqslant p < 4$.
   Without loss of generality we can assume that $0< T < 1$. Then
   \begin{equation*}
      \big \| \mathscr N \big( u \big)(t) \big \|_{\calX_T} \leqslant \|u_0\|_{\calX_T} + C T^{1-\nicefrac{1}{p}} \big \| u(t) \big \|_{\calX_T}
   \end{equation*}
   and
   \begin{multline*}
      \rule{20pt}{0pt}\big \| \mathscr N \big( u \big)(t) - \mathscr N \big( u \big)(t) \big \|_{\calX_T} \\
      \leqslant C T^{1-\nicefrac{1}{p}} \big \| u(t) - v(t)\big \|_{\calX_T} \Big( \big \| u(t) \big \|_{\calX_T} + \big \| v(t) \big \|_{\calX_T} \Big). \rule{20pt}{0pt}
   \end{multline*}
   For arbitrary $R>0$ we can choose  $T>0$ so small that the operator $\mathscr N$ is a contraction in the ball $B_R(0) \subset \calX_T$. The Banach fixed point theorem gives the existence of a local in time solution in the space $\calX_T$. From the construction this solution is unique in the ball $B_R(0)$. But the uniqueness is a local property and the above solution is also unique in the space $\mathcal X_{\hat T}$ for $\hat T > 0$ possibly smaller than $T$. Indeed, let us suppose that there is another solution in the space $\calX_T$ contained in some ball $B_{\hat R}(0)$ with $R < \hat R$. Taking $\hat T < T$ we construct the unique solution on the time interval $[0,\hat T]$ but in the ball $B_{\hat R}(0)$, which contradicts the uniqueness of the previous solution{.\hskip 10pt$\blacksquare$}
   \medskip

\begin{proposition}\label{AGHnieujemnosc}
   Let $\nicefrac{4}{3} \leqslant p < 4$ and $u_0(x) \in L^1(\R2) \cap L^p(\R2)$. If additionally $u_0(x)$ is nonnegative, then for all $0 < t <T$, $u(x,t) \geqslant 0$ and
   \begin{equation*}
      \int_{\R2} u(x,t) \dx = \int_{\R2} u_0(x) \dx.
   \end{equation*}
\end{proposition}

{\hskip -\parindent \sc Proof \hskip 7pt}
   The conservation of positivity for the solutions of the problem \eqref{AGHproblem} is a consequence of the analogous property for the system on bounded domains. Indeed, there has been proved in \cite[Proposition 1 and Proposition 2]{biler_92} that solutions with initial conditions $0 \leqslant u_0(x) \in L^1(\Omega) \cap L^p(\Omega)$ are nonnegative for $t>0$. The proof involved Stampacchia truncation method. The same approach permits to prove nonnegativity of the solution $u(x,t)$ either by approximating $u(x,t)$ on the sets $\big\{ |x| < N \big \}$ where $N\in \mathbb N$, by solutions $u(x,t) = u_k(x,t)$ of the problem posed on the sets $\big\{ |x| < N + k \big \}$ where $k\in \mathbb N$ or by a direct rewriting of the Stampacchia scheme (see also \cite[Lemma 1]{gajewski_groger} and \cite[Theorem I, Lemma 4.1]{gajewski}).

   The $L^1(\R2)$ norm of the solution is properly controlled, as it was justified in the proof of Proposition \ref{AGHlocex}, namely
      \begin{equation*}
      \begin{split}
         \sup_{0 \leqslant t < T}\|\mathcal N(u)(t)\|_1  \leqslant \|u_0\|_1 + C T^{1-\nicefrac{1}{p}} \big \| u(t) \big \|_{\calX_T}.
      \end{split}
      \end{equation*}
   Integrating the Duhamel formula \eqref{AGHOperLoc} over $\mathbb R^2$ we get
   \begin{equation*}
      \int_{\R2} u(x,t) \dx = \int_{\R2} \heat{t} u_0(x) \dx + \int_{\R2} \int_0^t \nabla \heat{t-\tau}\cdot \big(u \nabla \phi \big)(\tau) \dy{\tau} \dx.
   \end{equation*}
   But the last integral on the right hand side is equal to zero, since we can apply the  Fubini theorem, and the first one is equal to $\int_{\R2} u_0(x) \dx${.\hskip 10pt$\blacksquare$} \medskip

We obtain the global in time existence of the solutions using the following \textit{a~priori} estimate.
\begin{lemma}\label{AGHmaleniep}
   Let $u(x,t)$ be a solution of the problem \eqref{AGHproblem} with nonnegative initial condition $u_0(x)\in L^1(\mathbb R^2)$. Let $M = \int_{\mathbb R^2} u_0(x) \dx$. For every $p\in(1,\infty)$ such that $2^{k-1}<p<2^k$ for some $k\in\mathbb N$,  the inequality
\begin{equation*}
   \|u(\cdot,t)\|_p \leqslant C(k)\,M\, t^{-(1-\nicefrac{1}{p})}
\end{equation*}
holds with a constant $C(k)$ given explicitly. Moreover
\begin{equation*}
   \|u\|_\infty \leqslant C\,M\, t^{-1}
\end{equation*}
holds for some constant $C>0$.
\end{lemma}

We base the proof of the above lemma on the following interpolation inequality 
\begin{lemma}\label{AGHnesh}
Let $1\leqslant q,r \leqslant\infty,\ j,m\in\mathbb N,\ 0\leqslant j < m$ and $a\in[\nicefrac{j}{m},1]$. Moreover
   \begin{equation*}
      \frac{1}{p}=\frac{j}{n}+ a\Big(\frac{1}{r}-\frac{m}{n}\Big)+(1-a)\frac{1}{q},
   \end{equation*}
   where $n$ is the space dimension. Then 
   \begin{equation}\label{AGHnesh2}
      \sum_{|\alpha|=j} \|D^\alpha u\|_p\leqslant C \bigg( \sum_{|\alpha| = m} \|D^\alpha u\|_r\bigg)^a\,\|u\|_q^{1-a}.
   \end{equation}
\end{lemma}
This inequality was introduced by Gagliardo-Nirenberg-Sobolev and sometimes is referred to as the Nash inequality \cite{nash}.

{\hskip -\parindent \sc Proof of Lemma \ref{AGHmaleniep}\hskip 10pt}\\
   \textsc{Step 1. \hskip 5pt}
      Let $p=2$. Multiplying the first equation in \eqref{AGHmainintro} by $u$ and integrating by parts, we get 
      \begin{equation}\label{AGHrachunek}
      \begin{split}
         0= &\int_{\mathbb R^2} \big( uu_t-u \Delta u - u \nabla \cdot (u \nabla \phi )\big) \dx\\
         &=\frac{1}{2}\frac{\text d}{\text dt}\|u\|^2_2+\|\nabla u\|^2_2 +\frac{1}{2}\int_{\mathbb R^2}\nabla u^2\cdot\nabla\phi\dx\\
         &=\frac{1}{2}\frac{\text d}{\text dt}\|u\|^2_2+\|\nabla u\|^2_2 +\frac{1}{2}\int_{\mathbb R^2} u^3 \dx.
      \end{split}
      \end{equation}
      From Proposition \ref{AGHnieujemnosc} the solution $u(x,t) \geqslant 0$ in $\mathbb R^2$ for all $t>0$. Then we obtain the inequality
      \begin{equation}\label{AGH12}
         \frac{\text d}{\text dt}\|u\|^2_2+2\|\nabla u\|^2_2 \leqslant 0.
      \end{equation}
      Then taking $a=\nicefrac{1}{2},\ q=m=1,\ n=r=p=2,\ j=0$ in \eqref{AGHnesh2} we get
      \begin{equation}\label{AGHnesh5}
         \|u\|_2\leqslant C\|\nabla u\|^{\nicefrac{1}{2}}_2 \|u\|_1^{\nicefrac{1}{2}}
      \end{equation}
      which is equivalent to
      \begin{equation}\label{AGHnesh3}
         \frac{1}{CM_0^{2}}\|u\|_2^4\leqslant \|\nabla u\|^{2}_2
      \end{equation}
      where $M_0=\|u(0)\|_1$. Using \eqref{AGH12} and \eqref{AGHnesh3} we get
      \begin{equation*}
         \frac{\text d}{\text dt} \|u\|_2^2 + \frac{2}{CM_0^2}\|u\|_2^4 \leqslant 0.
      \end{equation*}
      Integrating the above inequality yields
      \begin{equation*}
         \|u(t)\|_2^2\leqslant \frac{1}{\frac{1}{\|u(0)\|_2^2} + \frac{2}{CM_0^2}t} \leqslant \frac{CM_0^2}{2}\:t^{-1},
      \end{equation*}
      which ends this part of the proof.\\
   \textsc{Step 2. \hskip 5pt}
      $p=2^k$. We prove it by recurrence. Let as assume that
      \begin{equation}\label{AGHzalind}
         \|u\|_{2^{k-1}}^{2^{k-1}}\leqslant a_{k-1}M_0^{2^{k-1}}t^{-(2^{k-1}-1)},
      \end{equation}
      where
      \begin{equation}
         a_k=\begin{cases}
                \nicefrac{C}{2} & \text{\quad for\quad } k=1,\\
                C2^{k-2}\big(a_{k-1}\big)^2 & \text{\quad for\quad } k \geqslant 2.
             \end{cases}
      \end{equation}
      In the inequality \eqref{AGHnesh5} we set $u^{2^{k-1}}$ instead of $u$ which yields
      \begin{equation}\label{AGHnesh7}
         \big(\|u\|_{2^k}^{2^{k}}\big)^2\leqslant C\|\nabla u^{2^{k-1}}\|_2^2 \big(\|u\|^{2^{k-1}}_{2^{k-1}}\big)^2.
      \end{equation}
      The same calculation as in \eqref{AGHrachunek} with the multiplication by $u^{2^{k}-1}$ gives
      \begin{equation}\label{AGHrach2}
         \frac{\text d}{\text dt}\|u\|_{2^k}^{2^k}+\frac{4(2^k-1)}{2^k}\|\nabla u^{2^{k-1}}\|_2^2\leqslant 0.
      \end{equation}
      Combining together inequalities \eqref{AGHzalind}, \eqref{AGHnesh7} and \eqref{AGHrach2} we get
      \begin{gather*}
         \|\nabla u^{2^{k-1}}\|_2^2 \geqslant  \frac{\big(\|u\|_{2^{k}}^{2^{k}}\big)^2}{C \big(\|u\|_{2^{k-1}}^{2^{k-1}}\big)^2},\\
         \frac{\text d}{\text dt}\|u\|_{2^k}^{2^k}+\frac{2^k-1}{C  2^{k-2}} \frac{\big(\|u\|_{2^{k}}^{2^{k}}\big)^2}{\big(\|u\|_{2^{k-1}}^{2^{k-1}}\big)^2}\leqslant 0,\\
         \frac{\frac{\text d}{\text dt}\|u\|_{2^k}^{2^k}}{\big(\|u\|_{2^{k}}^{2^{k}}\big)^2} \leqslant -\frac{2^k-1}{C M_0^{2^k} 2^{k-2}(a_{k-1})^2}\:t^{2^k-2} \hskip 3pt .
      \end{gather*}
      Let us denote $f(t)=\|u(t)\|_{2^k}^{2^k}$ and integrate the above inequality over the interval $(0,t)$. We obtain
      \begin{gather*}
         -\frac{1}{f(t)}+\frac{1}{f(0)}\leqslant -\frac{1}{C M_0^{2^k} 2^{k-2}\big(a_{k-1}\big)^2}\: t^{2^{k}-1},\\
         \frac{1}{f(t)}\geqslant \frac{1}{C  M_0^{2^k} 2^{k-2}\big(a_{k-1}\big)^2} t^{2^{k}-1}+\frac{1}{f(0)} \geqslant \frac{1}{C  M_0^{2^k} 2^{k-2}\big(a_{k-1}\big)^2}\: t^{2^{k}-1},\\
         f(t)\leqslant C2^{k-2}\big(a_{k-1} M_0^{2^k} \big)^2 t^{-(2^k-1)}
      \end{gather*}
      which gives
      \begin{equation}\label{AGHwynik}
         \|u\|_{2^k}^{2^k}\leqslant a_k  M_0^{2^k} t^{-(2^k-1)}.
      \end{equation}
   \textsc{Step 3. \hskip 5pt}
      General case. For each $p\in (1,\infty)$  there exists some $k\in\mathbb N$ such that $2^{k-1}<p<2^k$. Then it is sufficient to use inequality \eqref{AGHwynik} and very well known interpolation inequality
      \begin{equation*}
         \|u\|_p \leqslant \|u\|_{2^k}^{1-a}\|u\|_1^a
      \end{equation*}
      with $a=\frac{1-\nicefrac{1}{p}}{1-\nicefrac{1}{2^k}}$. Then
      \begin{equation*}
         \|u\|_p \leqslant M_0 a_k^{(\nicefrac{1}{p}-\nicefrac{1}{2^k})/(2^k-1)}t^{-(1-\nicefrac{1}{p})}.
      \end{equation*}
      At the end let us notice that
      \begin{equation*}
         a_k^{(\nicefrac{1}{p}-\nicefrac{1}{2^k})/(2^k-1)}\leqslant a_k^{\nicefrac{1}{p}-\nicefrac{1}{2^k}}\leqslant a_k^{\nicefrac{1}{2^k}}
      \end{equation*}
   \textsc{Step 4. \hskip 5pt}
      $p=\infty$. In the last step we use the property
      \begin{equation*}
         \|u\|_\infty = \limsup_{p\,\rightarrow\,\infty} \|u\|_p.
      \end{equation*}
      Then it is enough to show that $\lim_{k\rightarrow\infty} a_k^{\nicefrac{1}{2^k}}<\infty$. But $a_k=C^{v_k} 2^{w_k}$, where
      \begin{equation*}
         w_k=\begin{cases}
             -1 & \text{\quad for\quad } k=1,\\
             2w_{k-1}+(k-2)& \text{\quad for\quad } k \geqslant 2.
          \end{cases}
      \end{equation*}
      and $v_k=2^k-1$. It is not difficult to see that for all $k\in\mathbb N$,  $w_k < 0$. Then
      \begin{equation*}
         \lim_{k \rightarrow \infty} a_k^{\nicefrac{1}{2^k}} = \lim_{k \rightarrow \infty} \frac{C^{1-\nicefrac{1}{2^k}}}{2^{\nicefrac{|w_k|}{2^k}}}\leqslant C<\infty
      \end{equation*}
      which ends the proof {.\hspace{10pt}$\blacksquare$}
   \medskip

As a simple corollary from Proposition \ref{AGHlocex} and Lemma \ref{AGHmaleniep} we obtain the main result of this section, namely 
\begin{theorem}\label{AGHGlobExLp}
   Let $\nicefrac{4}{3} \leqslant p < 4$ and $u_0(x) \in L^1(\R2) \cap L^p(\R2)$. Then, there exists the unique solution of the integral problem \eqref{AGHcalkowe} in the space $C\big( [0,\infty):\,L^1(\R2) \big) \cap C\big( [0,\infty):\,L^p(\R2) \big)$.
\end{theorem}

\section{Radially symmetric self-similar solutions}\label{AGHsecradial}

As we mentioned in the introduction the system \eqref{AGHmainintro} is invariant under the scaling 
\begin{equation}\label{AGHscalling}
   \lambda^2 u(\lambda x,\lambda^2 t),\quad \phi(\lambda x,\lambda^2 t),
\end{equation}
where $0<\lambda\in\mathbb R$.

In this section we assume that $u(x,t),\ \phi(x,t)$ are self-similar, i.e.\ they are invariant under the scaling \eqref{AGHscalling}. The most important from our point of view are the solutions beginning their evolution from the initial conditions of the form $M \delta_0(x)$. Such an initial condition is obviously radial. Thus we assume from now on the radial symmetry of the solutions $u(x,t),\ \phi(x,t)$.\\
Following {\it the integrated density} method (see \cite{biler_95, biler_nadzieja_raczynski,olech_04}) we introduce the new variable
\begin{equation*}
   Q(r,t)=\int _{B_r(0)} u(x,t)\dx=2\pi \int _0^r s\,u(s,t)\dy{s},
\end{equation*}
where $B_r(0)$ is the closed ball of radius $r$ centered at the origin $(0,0)$. Integrating the system \eqref{AGHmainintro} over such balls and performing simple computations, we rewrite the problem \eqref{AGHmainintro}--\eqref{AGHinitintro} in the form
\begin{gather*}
   Q_t = Q_{rr} - \frac{1}{r}Q_r -\frac{1}{2\pi r} Q Q_r,\\
   Q(0,t) = 0, \quad \text{for} \quad t \geqslant 0\\
   \lim_{r \rightarrow \infty} Q(r,t) = \int_{\mathbb R^2} u_0(x)\dx = M, \quad \text{for} \quad t \geqslant 0.
\end{gather*}
Using self-similarity we can again change the variables: $Q(r,t) = 2\pi \xi (y)$ for $y=\nicefrac{r^2}{t}$. This way we transform our problem to the following one 
\begin{subequations}\label{AGHxi_intro}
   \begin{align}
      &\xi''(y) +\frac{1}{4}\xi'(y)-\frac{1}{2y}\xi(y)\xi'(y)=0,\label{AGHnaxi_eq}\\
      &\xi(0)=0, \quad \xi'(0)=a.\label{AGHnaxi_in}
   \end{align}
\end{subequations}

Formally, the boundary conditions implied by \eqref{AGHinitintro} are \begin{equation*}
   \xi(0)=0 \quad \text{and} \quad \lim_{y\rightarrow\infty}\xi(y)=M.
\end{equation*}
But it is more convenient to investigate equation \eqref{AGHnaxi_eq} with the boundary conditions \eqref{AGHnaxi_in} for some positive real $a$. As we will see, both problems are equivalent. The reader can find similar results considering the problem \eqref{AGHxi_intro} (existence of solutions, regularity, etc.) in \cite{biler_nadzieja_raczynski} but methods used there are different.

The main result we prove in this paper is the following
\begin{theorem}\label{AGHmain_th_rss}
   For all $M>0$, there exists a unique solution of \eqref{AGHnaxi_eq}, satisfying the initial conditions $\xi(0)=0$ and $\lim_{y\rightarrow\infty}\xi(y) = M$.
\end{theorem}

First we justify some \textit{a~priori} properties of the solutions of \eqref{AGHxi_intro}. Then we are able to construct a suitable space where the Banach fixed point theorem can be applied.
\begin{lemma}\label{AGHwlasnosci_xi}
   Let $\xi(y)$ be a solution of \eqref{AGHxi_intro} with $y\in [0,Y_0)$. Then
   \begin{itemize}
      \item[i)] $\xi'(0)>\nicefrac{1}{2}$ implies $\xi'(y)>\nicefrac{1}{2}$, for all $y\in[0,Y_0)$;
      \item[ii)] if $0<\xi'(0) = a<\nicefrac{1}{2}$ then $0< \xi'(y) < a$ and $-\nicefrac{1}{8} < \xi''(y) < 0$ for all $y\in[0,Y_0]$. Moreover, in such a case $0 < \xi(y) < a\,y$ on the whole existence interval;
      \item[iii)]  if $\xi(y)$ is global and $0<\xi'(0)<\nicefrac{1}{2}$, then $\lim_{y\rightarrow\infty} \xi(y)$ exists;
   \end{itemize}
\end{lemma}

{\hskip -\parindent \sc Proof \hskip 10pt}
   {\it (i)\quad}
   In order to compute $\xi''(0)$, we let $y$ tend to $0^+$ in the equation \eqref{AGHnaxi_eq}. Then
   \begin{equation*}
      \xi''(0) = \xi'\big(0\big) \bigg( \frac{1}{2} \,\xi'(0) -\frac{1}{4} \bigg) =\frac{a}{2} \bigg(a-\frac{1}{2}\bigg)
   \end{equation*}
   since
   \begin{equation*}
      \xi'(0)=\lim_{y\rightarrow 0^+}\frac{\xi(y)-\xi(0)}{y}=\lim_{y\rightarrow 0^+}\frac{\xi(y)}{y}.
   \end{equation*}
   In the case $a>\nicefrac{1}{2}$ we have $\xi''(y)>0$ for some right neighbourhood of $y=0$. Now let us define
   \begin{equation}\label{AGHdef_y_0}
      y_0=\sup \Big \{ y\in[0,Y_0):\ \xi'(s)>\tfrac{1}{2} \text{\quad for all\quad } s \in [0,y)\Big \},
   \end{equation}
   and suppose that we have $\xi'(y_0)=\nicefrac{1}{2}$ at a point $y_0 < Y_0$.

   It follows immediately that $\xi'(y)$ is decreasing on some subinterval of $[0,y_0)$.
   Since for $y\in[0,y_0)$ we have $\xi(y) > \tfrac{1}{2}y$, then by \eqref{AGHnaxi_eq}, the second derivative $\xi''(y)$ is positive on the interval $[0,y_0)$. It follows that $\xi'(y)$ cannot decrease on this interval, which contradicts our assumption. Thus, we have $\xi'(y)>\nicefrac{1}{2}$ for all $y\in[0,Y_0)$.

   {\it (ii)\quad}
   It is easy to check, that if $0< a <\nicefrac{1}{2}$ then $\xi'(y)< \nicefrac{1}{2}$ on $[0,Y_0)$. It is enough to note that $\xi''(y) < 0$, which implies that the first derivative $\xi'(y)$ is strictly decreasing, so $\xi'(y) < \xi'(0) = a$. It remains to prove that $\xi'(y)$ is positive.
   Integrating $\xi'(y)<\nicefrac{1}{2}$, we get that $\xi(y)<\nicefrac{y}{2}$ for $y\in [0,Y_0)$. Assume on the contrary, that there exists a point $\tilde y\in [0,Y_0)$, such that $\xi'(\tilde y)=0$. Then $\xi''(\tilde y)=0$. We have to consider two possibilities: {\it 1.}~~$\xi'(y) < 0$; {\it 2.}~~$\xi'(y) > 0$; both, for $y$ in a small, right neighbourhood of $\tilde y$, i.e.\ $[\tilde y,\tilde y+\varepsilon) =: \mathcal R_\varepsilon^+(\tilde y)$, with sufficiently small $\varepsilon >0$.
   \begin{itemize}
      \item[\it ad 1.] We know that $\xi'(y)$ is decreasing in $\mathcal R_\varepsilon^+(\tilde y)$. On the other hand,  since
      \begin{equation*}
         \xi''(y)=\xi'(y)\left(\frac{\xi(y)}{2y}-\frac{1}{4}\right)>0,
      \end{equation*}
      and $\xi(y) < \nicefrac{y}{2}$ we deduce that $\xi''(y) > 0$ for $y>\tilde y$. Then the function $\xi'(y)$ cannot decrease for any $y>\tilde y$;
      \item[\it ad 2.] Here $\xi'(y)$ is increasing in $\mathcal R^+_\varepsilon(\tilde y)$. But in the same manner as in
      the previous case, we see that $\xi''(y) < 0$.
   \end{itemize}
   A contradiction ends this part of the proof since integrating the above inequality and using the initial conditions we obtain the remaining part of the conclusion.

   {\it (iii)\quad} We have shown that $\xi'(y)>0$ and $\xi''(y)<0$. We have to consider two possibilities: either the limit of $\xi'(y)$ at $\infty$ is equal to $0$ or it is a positive constant less than $\nicefrac{1}{2}$. We will exclude the second case.

   First of all we rewrite \eqref{AGHnaxi_eq} as
   \begin{equation*}
      \frac{\xi''(y)}{\xi'(y)}=\left(\frac{\xi(y)}{2y}-\frac{1}{4}\right).
   \end{equation*}
   Using {\it (ii)} we show that the right hand side of the above equality is strictly negative, namely
   \begin{equation}\label{AGHujemne_xi}
      \frac{\xi''(y)}{\xi'(y)}\leqslant -\delta = -\frac{1}{2}\bigg( \frac{1}{2} - a\bigg)<0.
   \end{equation}
   We can write \eqref{AGHujemne_xi} in the form
   \begin{equation*}
      \left(\ln \xi'(y)\right)'\leqslant -\delta,
   \end{equation*}
   and after the integration from  $0$ to some $y>y_0$, we obtain
   \begin{equation*}
   \begin{split}
      \ln\xi'(y) & < \ln\xi'(0) - \delta y,\\
         \xi'(y) & < \ee^{\ln a}\ee^{- \delta y},\\
         0<\xi'(y) & < a \ee^{ - \delta y}.
   \end{split}
   \end{equation*}
   Let $y$ go to $\infty$. Then
   \begin{equation*}
      \lim_{y\rightarrow\infty}\xi'(y)\leqslant a \lim_{y\rightarrow\infty} \ee^{- \delta y}\leqslant 0.
   \end{equation*}
   Integrating again the inequality $\xi'(y) < a \ee^{- \delta y}$ on the interval $[0,y]$, we obtain
   \begin{equation}
   \begin{split}\label{AGHoszac_xi}
      \xi(y) - \xi(0) & < a \int _0^y \ee^{- \delta s} \dy{s},\\
      \xi(y) & < \frac{a}{\delta}\big(1 - \ee^{- \delta y} \big) = \frac{4 a}{1-2a} \big( 1 - \ee^{\nicefrac{1}{2}(a-\nicefrac{1}{2})y} \big),
   \end{split}
   \end{equation}
   which concludes the argument. \hspace{10pt}$\blacksquare$\medskip

To prove the existence of solutions of the problem \eqref{AGHxi_intro} we transform it into an integral one and apply a Banach fixed point theorem. Multiplying equation \eqref{AGHnaxi_eq} by $\ee^{\nicefrac{y}{4}}$, integrating and using the initial conditions \eqref{AGHinitintro}, we obtain
\begin{align*}
   \Big(\ee^{\nicefrac{y}{4}}\xi'(y)\Big)' =&\ \frac{\ee^{\nicefrac{y}{4}}}{2y}\xi(y)\xi'(y), \\
   \xi'(y) =&\ a\ee^{-\nicefrac{y}{4}}+\frac{1}{2}\ee^{-\nicefrac{y}{4}}\int_0^y \frac{\ee^{\nicefrac{s}{4}}}{s}\,\xi(s)\xi'(s)\dy{s}, \\
   \xi(y)=&\ 4a \Big(1- \ee^{-\nicefrac{y}{4}}\Big) +\frac{1}{2}\int_0^y \bigg(\ee^{-\nicefrac{t}{4}}\int_0^t \frac{\ee^{\nicefrac{s}{4}}}{s} \,\xi(s) \xi'(s) \dy{s} \bigg) \dy{t}.
\end{align*}
We look for fixed points of the operator
\begin{equation*}
   \mathcal H(\xi)(y)=4a \Big(1- \ee^{-\nicefrac{y}{4}}\Big) +\frac{1}{2}\int_0^y \bigg(\ee^{-\nicefrac{t}{4}}\int_0^t \frac{\ee^{\nicefrac{s}{4}}}{s}\,\xi(s)\xi'(s)\dy{s}\bigg)\dy{t},
\end{equation*}
in the space
\begin{equation*}
   \calT = C^1[0,Y_0]\cap \Big \{ \xi :\ \xi(0)=0,\ \xi'(0)=a,\ 0<\xi'(y)<\tfrac{1}{2}, \text{\ for\ }y\in[0,Y_0] \Big \},
\end{equation*}
endowed with the usual $C^1$ norm
\begin{equation*}
   \|\xi\|_{\calT} = \sup_{0 \leqslant y \leqslant Y_0} \xi(y)+\sup_{0 \leqslant y \leqslant Y_0} \xi'(y) = \|\xi\|_C + \|\xi\|_{C^1}\hskip10pt.
\end{equation*}

\begin{lemma}\label{AGHexistence_xi}
   There exists  $Y_0>0$ such that for a given $0 < a < \nicefrac{1}{2}$ the problem \eqref{AGHxi_intro} has a solution in the space $\calT$. Moreover this solution is unique.
\end{lemma}

{\hskip -\parindent \sc Proof \hskip 10pt}
   We widely use here the properties of $\xi(y)$, which follow from the choice of the space $\calT$ and Lemma \ref{AGHwlasnosci_xi}. Let us begin with
   \begin{equation*}
   \begin{split}
      \big \|\mathcal H(\xi) \big \|_C \leqslant &\ 4a \Big(1- \ee^{-\nicefrac{Y_0}{4}}\Big) +\frac{1}{2} \sup_{0\leqslant y \leqslant Y_0} \int_0^y \ee^{-\nicefrac{t}{4}} \int_0^t \frac{\ee^{\nicefrac{s}{4}}}{s} \,|\xi(s) \xi'(s)| \dy{s} \dy{t}\\
      \leqslant &\ 4a \Big(1- \ee^{-\nicefrac{Y_0}{4}}\Big) + \frac{1}{2} \sup_{0\leqslant y \leqslant Y_0} \int_0^y \ee^{-\nicefrac{t}{4}} \int_0^t \frac{\ee^{\nicefrac{s}{4}}}{s} \,|\xi(s)|\,\|\xi\|_{C^1} \dy{s} \dy{t}\\
      \leqslant &\ 4a \Big(1- \ee^{-\nicefrac{Y_0}{4}}\Big) + \frac{1}{2} \|\xi\|_{C^1}^2 \sup_{0\leqslant y \leqslant Y_0} \int_0^y \ee^{-\nicefrac{t}{4}} \int_0^t \ee^{\nicefrac{s}{4}} \dy{s} \dy{t}\\
      \leqslant &\ 4a \Big(1- \ee^{-\nicefrac{Y_0}{4}}\Big) + 2 \|\xi\|_{C^1}^2 \sup_{0\leqslant y \leqslant Y_0} \big( y+4\ee^{-\nicefrac{y}{4}}-4 \big)\\
      \leqslant &\ 4a \Big(1- \ee^{-\nicefrac{Y_0}{4}}\Big) + 2 Y_0 \|\xi\|_{C^1}^2\hskip10pt.
   \end{split}
   \end{equation*}
   Similarly, we have
   \begin{equation*}
   \begin{split}
      \|\mathcal H(\xi)\|_{C^1} = &\ \sup_{0\leqslant y \leqslant Y_0} \bigg| a\ee^{-\nicefrac{y}{4}} +\frac{1}{2}\frac{\text d}{\text dy}\bigg(\int_0^y \Big( \ee^{-\nicefrac{t}{4}} \int_0^t \frac{\ee^{\nicefrac{s}{4}}}{s} \,\xi(s) \xi'(s) \dy{s} \Big) \dy{t} \bigg) \Bigg|\\ 
      \leqslant &\ a +\frac{1}{2} \sup_{0\leqslant y \leqslant Y_0} \ee^{-\nicefrac{y}{4}} \int_0^y \frac{\ee^{\nicefrac{s}{4}}}{s} \,|\xi(s)|\, |\xi'(s)| \dy{s}\\
      \leqslant &\ a + \frac{1}{2} \|\xi\|_{C^1} \sup_{0\leqslant y \leqslant Y_0} \ee^{-\nicefrac{y}{4}} \int_0^y as\,\frac{\ee^{\nicefrac{s}{4}}}{s}\dy{s} =\\
      \leqslant &\ a + 2 a \|\xi\|_{C^1} \sup_{0\leqslant y \leqslant Y_0} \Big(1-\ee^{-\nicefrac{y}{4}} \Big) = a + 2a \Big(1 - \ee^{-\nicefrac{Y_0}{4}} \Big) \|\xi\|_{C^1}\hskip10pt ,
   \end{split}
   \end{equation*}
   which gives
   \begin{equation*}
      \|\mathcal H(\xi)\|_{\calT} \leqslant a + 2a \Big(1- \ee^{-\nicefrac{Y_0}{4}}\Big)\big(2 + \|\xi\|_{C^1}\big) + 2 Y_0 \|\xi\|_{C^1}.
   \end{equation*}
   To show the contraction property for the operator $\mathcal{H}$ we estimate
   \begin{equation*}
   \begin{split}
      \|\mathcal H(\xi)-\mathcal H(\eta)\|_{C} \leqslant &\ \frac{1}{2} \sup_{0\leqslant y \leqslant Y_0} \int_0^y \ee^{-\nicefrac{t}{4}} \int_0^t \frac{\ee^{\nicefrac{s}{4}}}{s} \,|\xi(s)\xi'(s)-\eta(s)\eta'(s)| \dy{s} \dy{t}\\
      \leqslant &\ \frac{1}{2} \sup_{0\leqslant y \leqslant Y_0} \int_0^y \ee^{-\nicefrac{t}{4}} \int_0^t \frac{\ee^{\nicefrac{s}{4}}}{s} \,|\xi(s)|\cdot |\xi'(s)-\eta'(s)| \dy{s} \dy{t}\\
       & + \frac{1}{2} \sup_{0\leqslant y \leqslant Y_0} \int_0^y \ee^{-\nicefrac{t}{4}} \int_0^t \frac{\ee^{\nicefrac{s}{4}}}{s} \,|\xi(s)-\eta(s)|\cdot |\eta'(s)| \dy{s} \dy{t} \\
      \leqslant &\ \frac{1}{2} \sup_{0\leqslant y \leqslant Y_0} \int_0^y \ee^{-\nicefrac{t}{4}} \int_0^t a\ee^{\nicefrac{s}{4}} \sup_{0 \leqslant s \leqslant Y_0}|\xi'(s)-\eta'(s)| \dy{s} \dy{t}\\
       & + \frac{1}{2} \sup_{0\leqslant y \leqslant Y_0} \int_0^y \ee^{-\nicefrac{t}{4}} \int_0^t \frac{1}{2}\,\ee^{\nicefrac{s}{4}} \sup_{0 \leqslant s \leqslant Y_0}|\xi'(s)-\eta'(s)| \dy{s} \dy{t} \\
      \leqslant & \ \| \xi-\eta \|_{C^1} \big( 2a + 1 \big) \sup_{0 \leqslant y \leqslant Y_0} \Big( y+4\ee^{-\nicefrac{y}{4}}-4 \Big)\\
      \leqslant &\ Y_0 \big( 2a + 1 \big) \|\xi -\eta \|_{C^1} \hskip 5pt .
   \end{split}
   \end{equation*}
   We estimate the term $\|\mathcal H(\xi)-\mathcal H(\eta)\|_{\mathcal C^1}$ in the same way, to get
   \begin{equation*}
      \|\mathcal H(\xi)-\mathcal H(\eta)\|_{\calT}\leqslant C(Y_0) \|\xi'-\eta'\|_{\mathcal T}.
   \end{equation*}
   Choosing a sufficiently small $Y_0$ we justify, as in Theorem ~\ref{AGHGlobalExistence}, the existence of the solution to the integral problem. To complete the proof we have to show that $\xi(y)$ constructed above solves the differential equation \eqref{AGHnaxi_eq}. For all $y\in (0,Y_0)$ our solution is a smooth function since the singularity in the coefficient $\nicefrac{1}{y}$ disapears. Then we can  consider the same equation but on some interval not containing zero. The uniqueness of the solutions ends the argument. \hspace{10pt}$\blacksquare$ \medskip

\begin{lemma}\label{AGHglob_ex_xi}
   The solution obtained in Lemma \ref{AGHexistence_xi} exists for all real $y\in[0,\infty)$.
\end{lemma}

{\hskip -\parindent \sc Proof \hskip 10pt}
   Let us assume that the solution $\overline\xi(y)$ is not a global one. Then, there exists a finite
   \begin{equation}\label{AGHsup_y}
      Y_{\text{max}}=\sup \{y:\ \overline\xi(y) \text{\ exists and is in the space\ } \mathcal T \}\quad.
   \end{equation}
   Using monotonicity and boundedness of the solution on the existence interval given by Lemma \ref{AGHwlasnosci_xi}, we can let $y$ tend to $Y_{\text{max}}$. Let us denote $\overline\xi(y_{\text{max}})$ by $\alpha$ and $\overline\xi'(y_{\text{max}})$ by $\beta$. Then
   \begin{align}
      &\xi''(y)+\frac{1}{4}\xi'(y)-\frac{1}{2y}\xi(y)\xi'(y)=0,\\
      &\xi(Y_{\text{max}})=\alpha, \quad \xi'(Y_{\text{max}})=\beta.
   \end{align}
   is an ordinary differential Cauchy problem which is not a singular one on the interval $[Y_{\text{max}},Y_0')$ for some $Y_0' > Y_{\text{max}}$. Then there exists the solution $\tilde\xi(y)$ of this problem, on the interval $[Y_{\text{max}},Y_0')$ (see, e.g.\ \cite{hartman}). Since the function
   \begin{equation*}
      \eta(y) = \begin{cases}
         \overline\xi(y) & \text{\quad for \quad} y\in[0,Y_0),\\
         \tilde\xi(y)& \text{\quad for \quad} y\in[Y_0,Y_0')
      \end{cases}
   \end{equation*}
   belongs to $\mathcal T$  on the interval $[0,Y_0')$, we get a contradiction with the maximality of $Y_{\text{max}}$. \hspace{10pt}$\blacksquare$ \medskip

Now we can prove the main result of this paper.

{\hskip -\parindent \sc Proof of Theorem \ref{AGHmain_th_rss} \hskip 10pt}
   Lemmas \ref{AGHexistence_xi} and \ref{AGHglob_ex_xi} assure the existence of the unique solution of the problem \eqref{AGHnaxi_eq}--\eqref{AGHnaxi_in} for $y\in [0,\infty)$ and $a\in[0,\nicefrac{1}{2})$. Moreover, from the inequality \eqref{AGHoszac_xi} we obtain the estimate from above for total mass of the solution, i.e.\  $\xi(+\infty)$.

   As in \cite[Theorem 21]{protter_weinberger}, we use the comparison principle for a second order ordinary differential equation, to show the estimates from below of our solution. In order to do this we have to regularize our problem and show that all the estimates are independent of the introduced regularization parameter. Indeed, let $\xi_\eps$ solve the equation 
   \begin{equation}\label{AGHzregularyzowane}
      F(\xi_\eps) := \xi_\eps''(y)+\frac{1}{4}\xi_\eps'(y)-\frac{1}{2y+\eps} \;\xi_\eps(y)\; \xi_\eps'(y) = 0,
   \end{equation}
   with the initial conditions $\xi_\eps(0)=0,\ \xi_\eps'(0)=\theta$ for some positive $\eps$ and $\theta$. Then for $h(y) = A\big( 1 - \ee^{-\mu y} \big)$, with $A=\tfrac{4a}{1-2a}$ 
   \begin{equation*}
   \begin{split}
      F \big( h(y) \big) = & -\mu^2 A \ee^{-\mu y} -A\mu\ee^{-\mu y} \Big( \frac{1}{2y+\eps} A \big( 1-\ee^{-\mu y}\big)-\frac{1}{4}\Big)\\
      = & -\frac{\mu A}{4(2y+\eps)} \ee^{-\mu y}\Big( (2y +\eps) (4\mu - 1) + 4A \big( 1-\ee^{-\mu y}\big) \Big),
   \end{split}
   \end{equation*}
   which is not positive for all $y\in\mathbb R^+$, uniformly with respect to $\eps$. Existence of the solutions to the problem \eqref{AGHzregularyzowane} is a standard fact (see, e.g., \cite{hartman}). It suffices to show that family $\{\xi_\eps\}_{0<\eps<1}$ is equicontinuous on every finite interval $[0,Y_0]$. Then, passing to the limit $\eps\rightarrow 0$ we conclude the proof.\hskip 10pt$\blacksquare$ \medskip

   Now we see explicitly the equivalence of the both formulations of our problem: with the initial conditions $\xi(0)=0,\ \xi'(0)=a$ and the boundary conditions $\xi(0)=0,\ \xi(\infty)=M$. Taking $\xi'(0) =a$ we solve our system with total mass $\tfrac{4a}{1-2a}$ and if we want to know solutions with mass M we have to take $\xi'(0)=\tfrac{M}{2M+4}$ in \eqref{AGHnaxi_in}.

\section{Asymptotics}\label{AGHasym}

Before proving the main result of this section we recall one more technical fact (for the proof see \cite[Lemma 6.1]{karch}.
\begin{lemma}\label{AGHtechn}
   Let $h(x)$ be a nonnegative function in $L^1(0,1)$ such that $\int_0^1 h(x)\,{\rm d}x < 1$. Moreover, we assume that $\alpha(t),\ \beta(t):\ (0,\infty) \rightarrow \mathbb R^+$ are bounded and such that
   \begin{equation}\label{AGHwithint}
      \alpha(t) \leqslant \beta(t) + \int_0^1 h(x) \alpha(xt)\,{\rm d}x,
   \end{equation}
   for all $t\in (0,\infty)$.
   Then $\lim_{t\rightarrow\infty} \beta(t) = 0$ implies $\lim_{t\rightarrow\infty} \alpha(t) = 0$.
\end{lemma}

\medskip

In the case of the problem \eqref{AGHproblem} we are able to prove analogue of \cite[Theorem 3]{karch2}, namely the following theorem.
\begin{theorem}\label{do_wniosku}
   Under the assumptions of Theorem \ref{AGHGlobalExistence}, let $u,\ v$ be two solutions of the problem \eqref{AGHproblem} with the initial conditions $u_0,\ v_0 \in B\calX$, respectively, obtained in Theorem \ref{AGHGlobalExistence}. Then, for sufficiently small $\eps>0$ the relations
   \begin{equation}\label{AGHin11}
      \lim_{t \rightarrow \infty} t^{\nicefrac{1}{4}} \big \| \heat{t}(u_0 - v_0) \big \|_{\nicefrac{4}{3}} = 0
   \end{equation}
   and
   \begin{equation}\label{AGHin22}
      \lim_{t \rightarrow \infty} t^{\nicefrac{1}{4}} \big \| u(\cdot,t) - v(\cdot,t) \big \|_{\nicefrac{4}{3}} = 0
   \end{equation}
   are equivalent.
\end{theorem}

{\hskip -\parindent \sc Proof \hskip 10pt}
   The integral representation of the solution \eqref{AGHcalkowe} yields
   \begin{align*}
      u(t) - v(t) &= \heat{t}(u_0 - v_0) + \int_0^t \heat{t-\tau} \big( \calB (u,u) - \calB (v,v) \big)(\tau) \dy{\tau}\\
      &= \heat{t}(u_0 - v_0) + \int_0^t \heat{t-\tau} \big( \calB(u,u - v) - \calB (v,u - v) \big)(\tau) \dy{\tau}.
   \end{align*}
   Then
   \begin{multline}\label{AGHgwiazdka}
      \big \| u(t) - v(t) \big \|_{\nicefrac{4}{3}} \leqslant \big \| \heat{t}(u_0 - v_0) \big \|_{\nicefrac{4}{3}} \\
      + C  \int_0^t (t-\tau)^{-\nicefrac{3}{4}} \Big( \big \| u(\tau) \big \|_{\nicefrac{4}{3}} + \big \| v(\tau) \big \|_{\nicefrac{4}{3}} \Big) \big \| (u - v)(\tau)\big \|_{\nicefrac{4}{3}} \dy{\tau}.
   \end{multline}
   According to assumptions in Theorem \ref{AGHGlobalExistence}
   \begin{equation*}
      \sup_{t>0} t^{\nicefrac{1}{4}} \big \| u(t) \big \|_{\nicefrac{4}{3}} \leqslant 2 \eps \text{\quad and \quad} \sup_{t>0} t^{\nicefrac{1}{4}} \big \| v(t) \big \|_{\nicefrac{4}{3}} \leqslant 2\,\eps.
   \end{equation*}
   Multiplying the inequality \eqref{AGHgwiazdka} by $t^{\nicefrac{1}{4}}$ and denoting
   \begin{equation*}
      \alpha(t) = t^{\nicefrac{1}{4}} \big \| u(t) - v(t) \big \|_{\nicefrac{4}{3}},\quad \beta(t) = t^{\nicefrac{1}{4}} \big \| \heat{t}(u_0 - v_0) \big \|_{\nicefrac{4}{3}},
   \end{equation*}
   we get
   \begin{align*}
      \alpha(t) \leqslant&\ \beta(t) + C \int_0^t (t - \tau)^{-\nicefrac{3}{4}} t^{\nicefrac{1}{4}} \tau^{-\nicefrac{1}{4}} \times \\
      &\times \Big( \tau^{\nicefrac{1}{4}} \big \| u(\tau) \big \|_{\nicefrac{4}{3}} + \tau^{\nicefrac{1}{4}} \big \| v(\tau) \big \|_{\nicefrac{4}{3}} \Big) \big \| (u - v)(\tau)\big \|_{\nicefrac{4}{3}} \dy{\tau}\\
      =&\ \beta(t) + 4\, \eps\, C \int_0^1 (1-w)^{-\nicefrac{3}{4}} w^{-\nicefrac{1}{2}} \alpha(wt) \dy{w}.
   \end{align*}
   Since $h(w) = 4\, \eps\, C (1-w)^{-\nicefrac{3}{4}} w^{-\nicefrac{1}{2}} \in L^1(0,1)$ and, if necessary, at the expense of decreasing of $\eps>0$, the quantity $4\, \eps\, C \int_0^1 (1-w)^{-\nicefrac{3}{4}} w^{-\nicefrac{1}{2}} \dy{w}$ is strictly less then $1$, we can apply Lemma \ref{AGHtechn} to obtain that the relation \eqref{AGHin11} implies \eqref{AGHin22}.\\
   To obtain the reverse implication we proceed in a similar way. First let us estimate 
   \begin{equation}\label{AGHrevers}
   \begin{split}
      t^{\nicefrac{1}{4}} &\big \| u(t) - v(t) - \heat{t}\big( u_0 - v_0 \big) \big \|_{\nicefrac{4}{3}} \\
      &\leqslant C \int_0^1 (1-w)^{-\nicefrac{3}{4}} w^{-\nicefrac{1}{2}} (wt)^{\nicefrac{1}{4}} \big \| u(wt) - v(wt) \big \|_{\nicefrac{4}{3}} \dy{w}.
   \end{split}
   \end{equation}
   Then we notice that
   \begin{gather*}
      \sup_{w\in[0,1]} (wt)^{\nicefrac{1}{4}} \big \| u(wt) - v(wt) \big \|_{\nicefrac{4}{3}}
      \leqslant \sup_{t > 0} t^{\nicefrac{1}{4}} \big \| u(t) - v(t) \big \|_{\nicefrac{4}{3}}\\
      \leqslant \sup_{t > 0} t^{\nicefrac{1}{4}} \big \| u(t) \big \|_{\nicefrac{4}{3}} + \sup_{t > 0} t^{\nicefrac{1}{4}} \big \| v(t) \big \|_{\nicefrac{4}{3}} < \infty,
   \end{gather*}
   and use the Lebesgue dominated convergence theorem. We obtain that the right hand side of the inequality \eqref{AGHrevers} converges to zero as $t$ tends to $\infty$. Then we deduce that the relation \eqref{AGHin22} implies \eqref{AGHin11}.
   {\hskip 10pt $\blacksquare$}
   \medskip

   We derive as a corollary of Theorem \ref{do_wniosku} a remark on local asymptotic stability of self--similar solutions or, in other words, the self--similar asymptotics of solutions to \eqref{AGHcalkowe} with small initial data.
   \begin{corollary}
      Suppose that $u_0$ satisfies the assumptions of Theorem \ref{AGHGlobalExistence}, $\int_{\mathcal R^2} u_0(x) \dx = M$ and $v_0(x) = U_M(x,1)$, where $U_M$ is the unique self--similar solution with $\int_{\mathcal R^2} U_M(x,1) \dx = M$. Then
      \begin{equation*}
         \lim_{t \rightarrow \infty} t^{\nicefrac{1}{4}} \big \| u(\cdot,t) - U_M(\cdot,t) \big \|_{\nicefrac{4}{3}} = 0
      \end{equation*}
      holds.
   \end{corollary}

   {\hskip -\parindent \sc Proof \hskip 10pt}
      Indeed, if $v_0(x)$ is a regular initial condition with $\int_{\mathcal R^2} v_0(x) \dx = M$, then
      \begin{equation*}
         \lim_{t \rightarrow \infty} t^{\nicefrac{1}{4}} \big \| \heat{t}(u_0 - v_0) \big \|_{\nicefrac{4}{3}} = 0,
      \end{equation*}
      cf.\ \cite{karch}, because the $L^p$--asymptotics of the heat semigroup is determined by the Gauss kernel multiplied by $M$, namely
      \begin{equation*}
         \lim_{t \rightarrow \infty} t^{1 - \nicefrac{1}{p}} \big \| \heat{t} v_0 - M \heat{t} \delta_0 \big \|_{p} = 0, \quad 1 \leqslant p \leqslant \infty.
      \end{equation*}
      This is the reason why we took the self--similar profile at the time $t=1$, not the singular data at $t=0$.
   {\hskip 10pt $\blacksquare$}
   \medskip

   We expect that the local result in the above corollary can be extended to the global asymptotic stability of self--similar solutions on all $L^p$--spaces, that is 
      \begin{equation}\label{AGHGwaizdkaWniosek}
         \lim_{t \rightarrow \infty} t^{1 - \nicefrac{1}{p}} \big \| u(\cdot,t) - U_M(\cdot,t) \big \|_{p} = 0
      \end{equation}
   (which is in fact implied by \eqref{AGHGwaizdkaWniosek} for $p=1$), for each $u_0 \in B\calX$ and $U_M$ with $M = \int_{\mathcal R^2} u_0(x) \dx$.\\
   The proof of this result seems to be beyond the scope of methods in that paper. Indeed, imitating the variational approach as was in \cite{blanchet_dolb_pert} for the problem with the gravitational interactions (in particular: energy and entropy functionals as well as entropy dissipation relation), one may first prove a priori estimates on $u$ and then a bound on $\int_{\mathcal R^2} u(x,t) \log \frac{u(x,t)}{U_M(x,t)} \dx$, leading to \eqref{AGHGwaizdkaWniosek}.



\begin{thebibliography}{10}

\bibitem{babskii}
V.~Babskii, Z.~Zhukov, and V.~Yudovich.
\newblock {\em {M}athematical {T}heory of {E}lectrophoresis -- {A}pplications
  to {M}ethods of {F}ractionation of {B}iopolymers}.
\newblock Naukova Dumka, Kiev, 1983.

\bibitem{biler_92}
P.~Biler.
\newblock {\em {E}\-xi\-sten\-ce and a\-sym\-pto\-tics of so\-lu\-tions for a
  pa\-ra\-bo\-lic-elliptic system with nonlinear no-flux boundary condition}.
\newblock {Nonlinear~Analysis}, 19(12):1121--1136, 1992.

\bibitem{biler_95}
P.~Biler.
\newblock {\em The {Cauchy} problem and self-similar solutions for a nonlinear
  parabolic equation}.
\newblock {Studia~Math.}, 114(2):181--205, 1995.

\bibitem{biler_dolbeault}
P.~Biler and J.~Dolbeault.
\newblock {\em {L}ong time behavior of solutions of {Nernst}-{Planck} and
  {Debye}-{Hueckel} drift-diffusion systems}.
\newblock {Ann.~Henri~Poincar\'e}, 1(3):461--472, 2000.

\bibitem{biler_hebisch_nadzieja}
P.~Biler, W.~Hebisch, and T.~Nadzieja.
\newblock {\em The {Debye} system: existence and large time behavior of solutions}.
\newblock {Nonlinear~Analysis}, 23(9):1189--1209, 1994.

\bibitem{biler_nadzieja}
P.~Biler and T.~Nadzieja.
\newblock {\em Existence and nonexistence of solutions for a model of gravitational
  interaction of particles. {I}}.
\newblock {Colloq.~Math.}, 66(2):319--334, 1994.

\bibitem{biler_nadzieja_raczynski}
P.~Biler, T.~Nadzieja, and A.~Raczy{\'n}ski.
\newblock {\em Nonlinear singular parabolic equations}.
\newblock In {Reaction diffusion systems (Trieste, 1995)}, volume 194 of
  {Lecture Notes in Pure and Appl. Math.}, pages 21--36. Dekker, New York,
  1998.

\bibitem{blanchet_dolb_pert}
A.~Blanchet, J.~Dolbeault, and B.~Perthame.
\newblock {\em Two-dimensional {Keller}-{Segel} model: optimal critical mass and
  qualitative properties of the solutions}.
\newblock {Electron.\ J.\ Diff.\ Eq.}, 44, 2006.

\bibitem{choi_lui}
Y.~S. Choi and R.~Lui.
\newblock {\em Multi-dimensional electrochemistry model}.
\newblock {Arch.~Rational~Mech.~Anal.}, 130(4):315--342, 1995.

\bibitem{debye_huckel}
P.~Debye and E.~Hueckel.
\newblock {\em {Z}ur {T}heorie der {E}lectrolyte {II}}.
\newblock {Phys.~Zeit.}, 24:305--325, 1923.

\bibitem{gajewski}
H.~Gajewski.
\newblock {\em On existence, uniqueness and asymptotic behavior of solutions of the
  basic equations for carrier transport in semiconductors}.
\newblock {Z.~Angew.~Math.~Mech.}, 65(2):101--108, 1985.

\bibitem{gajewski_groger}
H.~Gajewski and K.~Groeger.
\newblock {\em On the basic equations for carrier transport in semiconductors}.
\newblock {J.~Math.~Anal.~Appl.}, 113(1):12--35, 1986.

\bibitem{hartman}
P.~Hartman.
\newblock {\em {O}rdinary {D}ifferential {E}quations}.
\newblock S. M. Hartman, Baltimore, 1973.
\newblock Corrected reprint.

\bibitem{horstmann_1}
D.~Horstmann.
\newblock {\em From 1970 until present: the {K}eller-{S}egel model in chemotaxis and
  its consequences. {I}}.
\newblock {Jahresber.~Deutsch.~Math.-Verein.}, 105(3):103--165, 2001.

\bibitem{horstmann_2}
D.~Horstmann.
\newblock {\em From 1970 until present: the {K}eller-{S}egel model in chemotaxis and
  its consequences. {II}}.
\newblock {Jahresber.~Deutsch.~Math.-Verein.}, 106(2):51--69, 2001.

\bibitem{jungel}
A.~J{\"u}ngel.
\newblock {\em {Q}uasi-hydrodynamic {S}emiconductor {E}quations}.
\newblock Progress in Nonlinear Differential Equations and their Applications,
  41. Birkh{\"a}user Verlag, Basel, 2001.

\bibitem{karch}
G.~Karch.
\newblock {\em Scaling in non\-li\-near pa\-ra\-bo\-lic equa\-tions}.
\newblock {J.\ Math.\ Anal.\ Appl.}, 234(2):534--558, 1999.

\bibitem{karch2}
G.~Karch.
\newblock {\em Scaling in nonlinear parabolic equations: applications to {D}ebye
  system}.
\newblock In {Disordered and complex systems (London, 2000)}, volume 553 of
  {AIP~Conf.~Proc.}, pages 243--248. Amer.~Inst.~Phys., Melville,~NY, 2001.

\bibitem{lieb_loss}
E.~H. Lieb and M.~Loss.
\newblock {\em Analysis}.
\newblock Graduate Studies in Mathematics. American Mathematical Society,
  Providence, RI, second edition, 2001.

\bibitem{markowich_schmeiser}
P.~A. Markowich, C.~A. Ringhofer, and C.~Schmeiser.
\newblock {\em {S}emiconductor {E}quations}.
\newblock Springer-Verlag, Vienna, 1990.

\bibitem{nash}
J.~Nash.
\newblock {\em Continuity of solutions of parabolic and elliptic equations}.
\newblock {Amer.\ J.\ Math.}, 80:931--954, 1958.

\bibitem{olech_04}
M.~Olech.
\newblock {\em Nonuniqueness of steady states in annular domains for {S}treater
  equations}.
\newblock {{Appl.\ Math.\ (Warsaw)}}, 31(3):303--312, 2004.

\bibitem{protter_weinberger}
M.~H. Protter and H.~F. Weinberger.
\newblock {\em {M}aximum {P}rinciples in {D}ifferential {E}quations}.
\newblock Springer-Verlag, New York, 1984.
\newblock Corrected reprint of the 1967 original.

\bibitem{suzuki_book}
T.~Suzuki.
\newblock {\em {F}ree {E}nergy and {S}elf-interacting {P}articles}.
\newblock Progress in Nonlinear Differential Equations and Their Applications,
  62. Birkh{\"a}user, Boston, MA, 2004.

\end{thebibliography}
\end{document}